\documentclass[a4paper,11pt]{article}

\usepackage[utf8]{inputenc}
\usepackage[english]{babel}
\usepackage[bitstream-charter]{mathdesign}
\usepackage[T1]{fontenc}
\usepackage{amsfonts}
\usepackage{geometry}
\usepackage{amsmath, amsthm,graphicx}
\usepackage{hyperref}
\usepackage{authblk}

\theoremstyle{plain}                   
\newtheorem{theorem}{Theorem}[section]  
\newtheorem{proposition}[theorem]{Proposition}   
        
\newtheorem{lemma}[theorem]{Lemma}           
\theoremstyle{definition}               
\newtheorem{definition}[theorem]{Definition}

\newtheorem*{Acknowledgement}{Acknowledgements}
\newtheorem*{Notation}{Notation}                    
\newtheorem{remark}[theorem]{Remark}

\newtheorem*{conjecture}{Conjecture}
\newtheorem*{example n=2}{Case n=2}

\date{}

\begin{document}

\title{A note on the Waring Ranks of Reducible Cubic Forms}

\author{Emanuele Ventura \thanks{Aalto University, Department of Mathematics and Systems Analysis, \textit{Email}: \texttt{emanuele.ventura@aalto.fi}}}

\maketitle

\begin{abstract}

Let $W_3(n)$ be the set of Waring ranks of reducible cubic forms in $n+1$ variables. We prove that  $W_3(n)\subseteq \lbrace 1,\ldots, 2n+1\rbrace$.  

\end{abstract}

\section{Introduction}

\indent

Let $K$ be an algebraically closed field of characteristic zero, let $V$ be a $(n+1)$-dimensional $K$-vector space and $F\in S^d V$, namely a homogeneous polynomial of degree $d$ in $n+1$ indeterminates. The {\bf Waring problem for polynomials} asks for the least value $s$ such that there exist linear forms $L_1, \ldots, L_s,$ for which $F$ can be written as a sum 
\begin{equation}
F=L_1^d+\ldots+L_s^d. 
\end{equation}

\noindent
This value $s$ is called the \textit{Waring rank}, or simply the \textit{rank}, of the form $F$, and here it will be denoted by $rk(F)$. The Waring problem for a \textit{general form} $F$ of degree $d$ was solved by Alexander and Hirschowitz, in their celebrated paper \cite {AH}. 

\begin{theorem}[{\bf Alexander-Hirschowitz \cite{AH}}]\label{alexander-hirschowitz}

A general form $F$ of degree $d$ in $n+1$ variables is the sum of $\lceil \frac{1}{n+1}\binom{n+d}{d}\rceil$ powers of linear forms, unless

\begin{description}

\item $d=2$, $s=n+1$ instead of $\lceil \frac{n+2}{2}\rceil$;

\item $d=3$, $n=4$ and $s=8$ instead of $7$;

\item $d=4$, $n=2,3,4$ and $s=6,10,15$ instead of $5,9,14$ respectively.

\end{description}

\end{theorem}

\begin{remark}

The assumption on the characteristic is not necessary, see \cite{IK} for more details. 

\end{remark}

\indent

The Waring problem in the case of a given homogeneous polynomial is far from being solved. A major development in this direction is made in \cite{CCG} where the rank of any monomial and the rank of any sum of pairwise coprime monomials are computed. \\
The present paper concerns with the Waring rank of reducible cubic forms. The main result of this work is the following theorem. 

\begin{theorem}\label{theo1}
Let $W_3(n)$ be the set of ranks of reducible cubic forms in $n+1$ variables, then 
$$
W_3(n)\subseteq \lbrace 1,\ldots,2n+1\rbrace.
$$
\end{theorem}

\section{The Apolarity}

In this section, we recall basic definitions and facts; see \cite{IK} and \cite{RS} for details. 
\indent

Let $K$ be an algebraically closed field of characteristic zero, $S=\bigoplus_{i\geq 0} S_i=K[x_0, \ldots, x_n]$ and $T=\bigoplus_{i\geq 0} T_i=K[\partial_0, \ldots, \partial_n]$ be the \textit{dual ring} of $S$ (i.e. the ring of differential operators over $K$). $T$ is an $S$-module acting on $S$ by differentiation
$$
\partial^\alpha (x^\beta)=\left\{\begin{array}{rl}
\alpha ! \binom{\beta}{\alpha} x^{\beta-\alpha} &\mbox{ if } \beta \geq \alpha\\
0 &\mbox{ otherwise. } \\
\end{array}\right.
$$

\noindent
where $\alpha$ and $\beta$ are multi-indices. The action of $T$ on $S$ is classically called {\bf apolarity}. Note that $S$ can also act on $T$ with a (dual) differentiation, defined by
$$
x^\beta (\partial^\alpha)=\beta ! \binom{\alpha}{\beta} \partial^{\alpha-\beta},
$$
if $\alpha \geq \beta$ and $0$ otherwise.\\

\noindent
In this way, we have a non-degenerate pairing between the forms of degree $d$ and the homogeneous differential operators of order $d$. Let us recall some basic definitions. 

\begin{definition}

Let $F \in S$ be a form and $D\in T$ be a homogeneous differential operator. Then $D$ is \textit{apolar} to $F$ if $D(F)=0$.

\end{definition}

\begin{definition}

For any $F\in S^dV$, we define the ideal $F^\perp=\lbrace D\in T | D(F)=0 \rbrace \subset T$, called the \textit{principal system} of $F$. If $F\in S^dV$, for every homogeneous operator $D\in T$ of degree $\geq d+1$, we have $D(F)=0$, or equivalently $D\in F^\perp$. The principal system of $F$ is a \textit{Gorenstein ideal}.

\end{definition}

\begin{definition}

Given a homogeneous ideal $I\subset T$, the Hilbert function $\mbox{HF}$ of $T/I$ is defined as 
$$
\mbox{HF}(T/I,i)=\dim_K T_i -\dim_K I_i.
$$
\noindent The first difference function $\Delta \mbox{HF}$ of the Hilbert function of $T/I$ is defined as 
$$
\Delta \mbox{HF}(T/I,i)=\mbox{HF}(T/I,i)-\mbox{HF}(T/I,i-1),
$$
\noindent where $\mbox{HF}(T/I,-1)$ is set to be zero. 
\end{definition}

\noindent
Now, we recall the key result of this section.

\begin{lemma}[Apolarity lemma]\label{ApolarityLemma}

A form $F\in S^dV$ can be written as

\begin{equation}
F=\sum_{i=1}^s L_i^d,
\end{equation}

\noindent
where $L_i$ are linear forms pairwise linearly indipendent, if and only if there exists an ideal $I \subset F^\perp$ such that $I$ is the ideal of a set of $s$ distinct points in $\mathbb P^n$, where these $s$ points are the corresponding points of the linear forms $L_i$ in the dual space $\mathbb P^{n*}$.
\end{lemma}

\indent

For a proof of apolarity lemma \ref{ApolarityLemma} see  for instance \cite{IK}. We will refer to the $s$ points of this lemma as \textit{decomposition points}.

\section{Classification of Ranks of Reducible Cubic Forms in $\mathbb P^n$}

\indent

In this section we give the classification of the ranks of reducible cubic forms. Since the rank is invariant under projective transformations, we only need to check the projective equivalence classes of cubic forms. Let $W_3(n)$ be the set of values of ranks of reducible cubic forms in $n+1$ variables, namely forms of type $F=LQ$, where $L,Q\in S$ are linear and quadratic forms respectively. In order to give a classification, note that $W_3(n-1)\subset W_3(n)$. Indeed, every form in $n$ indeterminates is also a form in the ring of polynomials in $n+1$ indeterminates and the ranks as polynomial in $n$ variables and as polynomial in $n+1$ variables are equal. The subset $W_3(n-1)\subset W_3(n)$ is the set of the ranks of reducible cones in $n+1$ variables. The forms $F=LQ$ which are not cones (up to projective equivalence) are the following.

\begin{itemize}

\item (Type $A$) $Q$ is not a cone and $L$ is not tangent to $Q$.

\begin{center}
\includegraphics[scale=0.1]{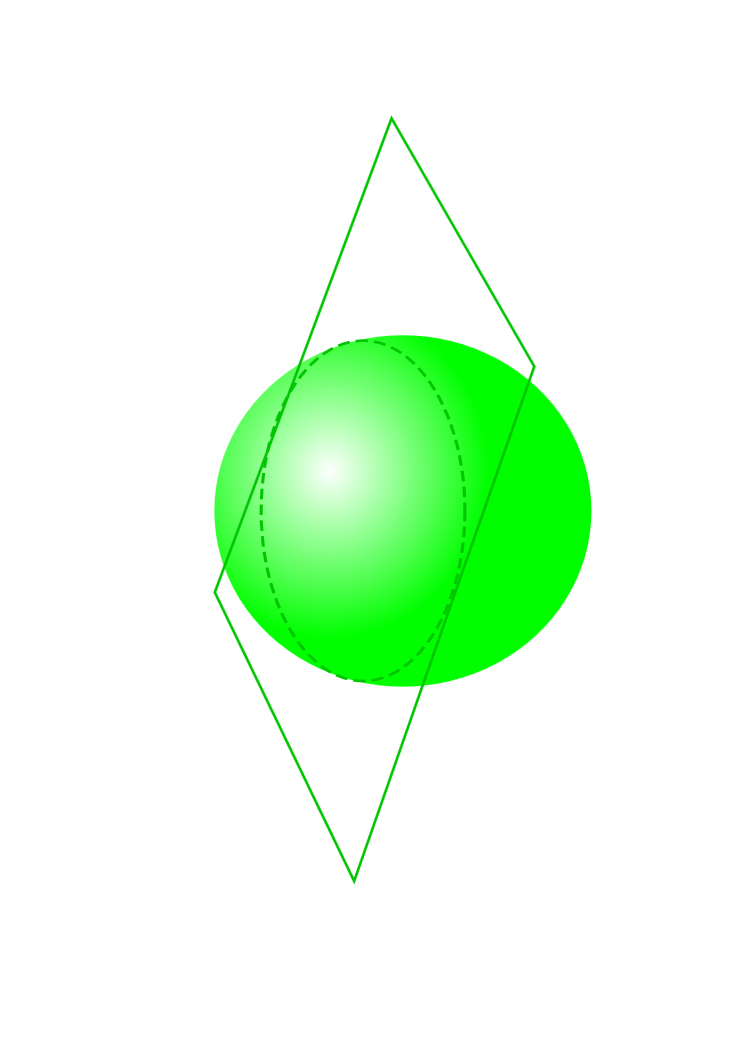}
\end{center} 
 
\item (Type $B$) $Q$ is a cone and $L$ does not pass through any vertex of $Q$.

\begin{center}
\includegraphics[scale=0.08]{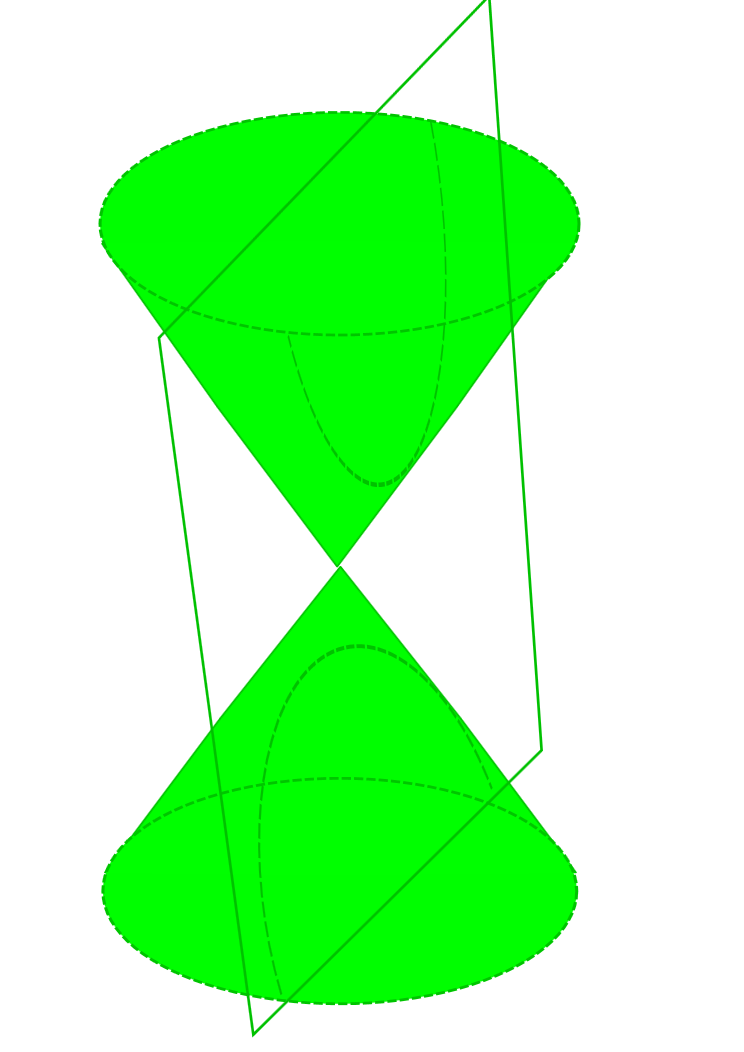}
\end{center} 

\item (Type $C$) $Q$ is not a cone and $L$ is tangent to $Q$.

\begin{center}
\includegraphics[scale=0.1]{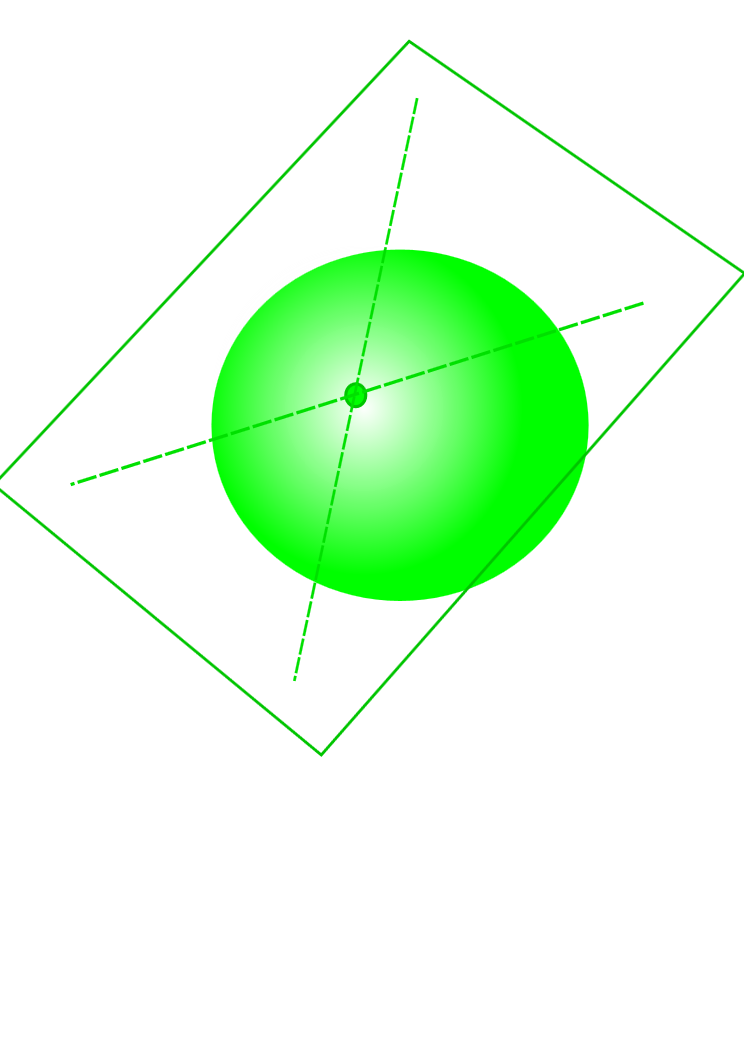}
\end{center} 

\end{itemize}

\newpage
We have the following result. 
\begin{theorem}\label{Theo2}
The ranks of reducible cubic forms $A$, $B$ and $C$ in $n+1$ variables are the following.
\begin{center}
    \begin{tabular}{ | l | l |}
    \hline
    Type & Rank  \\ \hline
     $A$ & $= 2n$   \\ \hline
     $B$ & $= 2n$   \\ \hline
     $C$ & $\geq 2n,\leq 2n+1$   \\ \hline
    \end{tabular}
\end{center}
\end{theorem}

The ranks of cubic forms of type $A$ and $B$ are given by [\cite{LT}, Proposition $7.2$]. B. Segre proved that the cubic surface in $\mathbb P^3$ of type $C$ has rank $7$ \cite{Segre}.  

\subsection{Type $C$}

Cubic forms of type $C$ are projectively equivalent to the cubic form

\begin{equation}\label{cubicform}
F=x_0(x_0x_1+x_2x_3+x_4^2+\ldots +x_n^2).
\end{equation}

\begin{Notation}

We denote by $\int G dx_i$ a suitable choice of a primitive of $G$ (that will be specified any time it is needed), namely a form $H$ such that $\partial_i H=G$, where $\partial_i$ denotes the usual partial derivative with respect to the variable $x_i$. 

\end{Notation}

First, note that if $n=2$, we have this proposition. 

\begin{proposition}\label{prop6}

The cubic form $F=x_0(x_0x_1+x_2^2)$ has rank $\leq 5$.
\begin{proof}
Consider the coordinate system given by the following linear transformation.
\begin{equation}\label{coordinate1}
\left\{\begin{array}{rl}
x_0=y_1\\
x_1=\frac{1}{3}y_1+y_3\\
x_2=y_2\\
\end{array}\right.
\end{equation}

\noindent
By this, we have $F=\frac{1}{3}y_1^3+y_1^2y_3+y_1y_2^2$. Let $K_1=\int \partial_2 F dy_2$ be the primitive of $\partial_2F$ given by $K_1=\frac{1}{6}[(y_1+y_2)^3+(y_1-y_2)^3]=\frac{1}{3}y_1^3+y_1y_2^2$. Thus $F=K_1+y_1^2y_3$, where $y_1^2y_3=\frac{1}{6}[(y_1+y_3)^3-(y_1-y_3)^3-2y_3^3]$. Then $rk(F)\leq 5$, which proves the statement.

\end{proof}

\end{proposition}

\noindent
It is straightforward to generalize this fact as follows.

\begin{proposition}\label{prop7}

The cubic form $F=x_0(x_0x_1+x_2x_3+x_4^2+\ldots +x_n^2)$ has rank $\leq 2n+1$.

\begin{proof}

We prove it by induction on $n$. The proposition holds for $n=2$ by Proposition \ref{prop6}.
Let us suppose the proposition true for all $i\leq n-1$ and prove the case $i=n$. Introduce the coordinate system given by the following linear transformation.
\begin{equation}\label{coordinate3}
\left\{\begin{array}{rll}
x_0=y_1\\
x_1=y_3\\
x_2=y_0+y_2\\
x_3=y_0-y_2\\
x_4=y_4\\
\vdots \ \ \ \ \  \\
x_n=y_n\\
\end{array}\right.
\end{equation}

\noindent
Then, the cubic becomes $F=y_0^2y_1-y_1y_2^2+y_1^2y_3+y_1y_4^2+\ldots +y_1y_n^2$.
Setting $G=\int \partial_0 F dy_0=y_0^2y_1+\frac{1}{3} y_1^3$, we take $F=G - \frac{1}{3} y_1^3-y_1y_2^2+y_1^2y_3+y_1y_4^2+\ldots+y_1y_n^2$. We have that $rk(G)=2$. Let $H=-\frac{1}{3}y_1^3-y_1y_2^2+y_1^2y_3+y_1y_4^2+\ldots +y_1y_n^2$. Since $H$ is a cubic form in $\mathbb P^{n-1}$ decomposed into a smooth quadric $Q$ and a tangent space $L$ to a point of $Q$ (and hence it is of type $C$), by inductive assumption $rk (H)\leq 2(n-1)+1$. Thus $rk (F)\leq rk (G) + rk (H)\leq 2+2(n-1)+1=2n+1$. Repeating the argument, one obtains a decomposition for $F$. 

\end{proof}

\end{proposition}

\begin{remark}

By [\cite{LT}, Theorem $1.3$], the rank of the cubic forms of type $C$ is $\geq 2n$. 

\end{remark}

\begin{remark}

The ranks of the reducible cubic forms are quite different from the generic rank of cubic forms given by the Alexander-Hirschowitz Theorem \ref{alexander-hirschowitz}: for sufficiently large values of $n$, the ranks of reducible cubics are smaller than the rank of the generic cubic.

\end{remark}

\section{Proof of Theorem \ref{theo1}}

\begin{proof}

We prove it by induction on $n$. If $n=1$, it is well known that cubic forms (actually, forms of any degree) in $2$ variables have rank at most their degree; in this case the set of ranks is exactly $W_3(1)$. Suppose that the statement holds for $i\leq n-1$ and we want to show it for $i=n$. Consider $ W_3(n)\setminus W_3(n-1)$; applying Theorem $\ref{Theo2}$, there exist forms of ranks $2n$ and of rank at most $2n+1$. By induction, $W_3(n-1)\subseteq \lbrace 1,\ldots,2n-1\rbrace$, and so $W_3(n)\subseteq \lbrace 1,\ldots,2n+1\rbrace$. 

\end{proof}

Motivated by the result of Segre \cite{Segre}, we state the following

\begin{conjecture}\label{conj on cubic form of type C}

The Waring rank of the reducible cubic forms of type $C$ in $n+1$ variables is $2n+1$.

\end{conjecture}

\begin{remark}\label{remark on HF}

The conjecture above states that $F=y_0^2y_1-y_1y_2^2+y_1^2y_3+y_1y_4^2+\ldots +y_1y_n^2$ has rank $\geq 2n+1$. The ideal $F^\perp$ is minimally generated by  $\partial_i\partial_3$ (for $i\neq 1$), $\partial_1\partial_3-\partial_i^2$ (for $i\neq 1,2,3$), $\partial_1\partial_3+\partial_2^2$, $\partial_i\partial_j$ (for $i,j\neq 1,3$), $\partial_i^3$ (for $i \neq 3$), $\partial_1^2\partial_i$ (for $i\neq 3$). \\
\indent The degree of a zero-dimensional scheme can be computed using Hilbert functions. Let $\mathbb X$ be a set of decomposition points of $F$ and set $I=I(\mathbb X)\subset F^\perp$. Let us suppose that $\mathbb X$ has no points on $\lbrace \partial_3=0\rbrace$. In this case, $\partial_3$ is not a zero-divisor in $T/I$, which is crucial here. Then the degree of $\mathbb X$ is given by

$$
\deg \mathbb X=\sum_{i\geq 0} \Delta \mbox{HF}(T/I,i)=\sum_{i\geq 0} \mbox{HF}(T/(I+\langle \partial_3\rangle),i)\geq \sum_{i\geq 0} \mbox{HF}(T/(F^\perp+\langle \partial_3\rangle),i)=2n+1,
$$

\noindent
where the Hilbert function $\mbox{HF}$ of $F^\perp+\langle\partial_3\rangle$ is the sequence $(1,n,n,0,-\cdots)$.\\
The case when $\mathbb X$ has points on $\lbrace \partial_3=0\rbrace$ requires a more careful analysis which we show for $n=2$. 
\end{remark}

We propose a technique based on apolarity and Hilbert functions that might be generalized to higher dimensions. We will show it dealing with the known case $n=2$. 

\begin{example n=2}
Let us denote $T=\mathbb C[\partial_1,\partial_2,\partial_3]$. In this case, we have $F=y_1(y_1y_3+y_2^2)$. The principal system of $F$ is the ideal 
$F^\perp=\langle \partial_1\partial_3-\partial_2^2,\partial_2\partial_3,\partial_3^2,\partial_1^3,\partial_1^2\partial_2,\partial_2^3\rangle$. Let $\mathbb X$ be a set of decomposition points of $F$ and let us set $I=I(\mathbb X)$.\\ 
\indent If $\mathbb X$ has no points on $\lbrace\partial_3=0\rbrace$ then $\partial_3$ is not a zero-divisor of $T/I$. Then  

$$\deg \mathbb X=\sum_{i\geq 0} \Delta \mbox{HF}(T/I,i)=\sum_{i\geq 0} \mbox{HF}(T/(I+\langle\partial_3\rangle),i)\geq \sum_{i\geq 0} \mbox{HF}(T/(F^\perp+\langle\partial_3\rangle),i)=5.$$
\noindent
Indeed, $I+\langle\partial_3\rangle\subset F^\perp+\langle\partial_3\rangle=\langle \partial_3,\partial_1\partial_2,\partial_2^2,\partial_1^3,\partial_2^3 \rangle$ and the Hilbert function of $T/(F^\perp+\langle\partial_3\rangle)$ is the sequence $(1,2,2,0,-\cdots)$, as in Remark \ref{remark on HF} above.\\
\indent Let us assume that $\mathbb X$ has some point on $\lbrace\partial_3=0\rbrace$. If $\dim I_2\leq 1$ then the Hilbert function of $T/I$ is the sequence $(1,3,m\geq 5,\ldots)$ and hence again $\deg \mathbb X\geq 5$. So let us assume $\dim I_2\geq 2$. Note that $I_2\subset F^\perp_2=\langle \partial_1\partial_3-\partial_2^2,\partial_2\partial_3,\partial_3^2\rangle$. 
There exists a two-dimensional subspace of conics  $L\subset I_2\subset F^\perp_2$. Either this space $L$ is the pencil $a\partial_3^2+b\partial_2\partial_3$, and the base locus of this pencil is $\lbrace\partial_3=0\rbrace$, or $I_2$ contains some irreducible conic of equation $\partial_1\partial_3-\partial_2^2+a\partial_3^2+b\partial_2\partial_3$, whose only common intersection with $\lbrace\partial_3=0\rbrace$ is the point $(1:0:0)$. The first case is not possible, since otherwise $\mathbb X\subset \lbrace\partial_3=0\rbrace$, namely $\partial_3 F=0$, which is false. Hence we have $\mathbb X\cap\lbrace\partial_3=0\rbrace=\{(1:0:0)\}.$ This implies that $\mathbb X\cap\lbrace\partial_3=0\rbrace\subset \mathbb X\cap\lbrace\partial_2=0\rbrace$. Then $\partial_3$ does not vanish at any point of $\mathbb X\cap\lbrace\partial_2\not=0\rbrace=\mathbb X'$. Note that $\deg \mathbb X'\leq \deg \mathbb X-1$ because the point $(1:0:0)$ does not belong to $\mathbb X'$. Setting $J=(I:\partial_2)$ the ideal of $\mathbb X'$, we have that $\partial_3$ is not a zero-divisor of $T/J$, so we can compute 
$$\deg \mathbb X'=\sum_{i\geq 0} \mbox{HF}(T/(J+\langle\partial_3\rangle),i)\geq \sum_{i\geq 0} \mbox{HF}(T/((F^\perp:\partial_2)+\langle\partial_3\rangle),i)\geq 4,$$
\noindent
since $(F^\perp:\partial_2)+\langle\partial_3\rangle=\langle \partial_3,\partial_1^2,\partial_2^2\rangle$. Finally $\deg \mathbb X\geq\deg \mathbb X'+1\geq  5$, which says that the rank of $F$ is at least $5$ using the apolarity lemma \ref{ApolarityLemma}. 

\end{example n=2}

\begin{Acknowledgement}

This paper is part of author's Master thesis at University of Catania. I would like to thank Riccardo Re for his support and encouragement. I also would like to thank Enrico Carlini and Zach Teitler for insightful discussions.

\end{Acknowledgement}



\begin{thebibliography}{99} 

\bibitem {AH} J. Alexander and A. Hirschowitz, \textit{Polynomial interpolation in several variables}, J. Alg. Geom. 4 (1995), n.2, 201-222.

\bibitem {CCG} E. Carlini, M.V. Catalisano, and A.V. Geramita, \textit{The solution to the Waring problem for monomials and the sum of coprime monomials}, J. of Algebra 370 (2012), 5-14.

\bibitem {IK} A. Iarrobino and V. Kanev, \textit{Power sums, Gorenstein Algebras, and Determinantal Loci}, volume 1721 of Lecture Notes in Mathematics. Springer-Verlag, Berlin (1999).

\bibitem {LT} J.M. Landsberg, Z. Teitler, \textit{On the ranks and border ranks of symmetric tensors}, Found Comput. Math. 10, (2010), no 3, 339–366.

\bibitem {RS} K. Ranestad and F.O. Schreyer, \textit{Varieties of Sums of Powers}, J. Reine Angew. Math. 525 (2000), 147–181.

\bibitem{Segre} B. Segre, \textit{The non-singular cubic surfaces}, Oxford University Press, Oxford, 1942. 

\end{thebibliography}
\end{document}